\documentclass{article}

\usepackage[utf8]{inputenc} % allow utf-8 input
\usepackage[T1]{fontenc}    % use 8-bit T1 fonts
\usepackage{url}            % simple URL typesetting

\usepackage{amsmath,amsthm}
\usepackage{amsmath}
\usepackage{amssymb}
\usepackage{amsfonts}
\usepackage{amsthm}
\usepackage{amsbsy}
\usepackage{amsgen}
\usepackage{amscd}
\usepackage{amsopn}
\usepackage{amstext}
\usepackage{amsxtra}
\usepackage{mathrsfs}
\usepackage{enumitem}
\usepackage{graphicx}
\usepackage{verbatim}
\usepackage{epstopdf}
\usepackage{float}
\usepackage[all,cmtip]{xy}
\usepackage{accents}
\usepackage{sseq}
\usepackage{url}
\usepackage{makeidx}
\usepackage{wrapfig}
\usepackage{tikz-cd}
\usepackage{framed,enumitem} 
\usepackage{syllogism}
\usepackage{tikz}
\usetikzlibrary{decorations.markings,arrows}

\usetikzlibrary{matrix}
\usepackage[T1]{fontenc}
\usepackage{titlesec, blindtext, color}
\definecolor{gray75}{gray}{0.75}
\newcommand{\hsp}{\hspace{20pt}}
\titleformat{\chapter}[hang]{\Huge\bfseries}{\thechapter\hsp\textcolor{gray75}{|}\hsp}{0pt}{\Huge\bfseries}

\usepackage[margin=0.6in]{geometry}

\theoremstyle{definition}

\theoremstyle{remark}

\newtheorem*{ex*}{Ejercicio}

\usepackage[square,sort,comma,numbers]{natbib}   % omit 'round' option if you prefer square brackets

\renewcommand{\contentsname}{Contenidos}

\usepackage{setspace}
\usepackage{afterpage}
\usepackage{sectsty}
\sectionfont{\large}

\AtBeginEnvironment{quote}{\small}

\usepackage{epigraph}
\usepackage{hyperref}       % hyperlinks

\begin{document}
\title{On a form of intrinsic optimism in Set Theory}
\author{M. Muñoz Pérez\footnote{\texttt{mupemiguel99@gmail.com}. Departamento de L\'ogica, Historia y Filosof\'ia de la Ciencia,
 UNED, Madrid, Spain.}}
\date{\today}

\maketitle

%1.25
\setstretch{1.25}
\setcounter{tocdepth}{5}
\renewcommand{\contentsname}{Contents}

\section*{Introduction}

In this paper we will try to provide a solid form of \textit{intrinsic set theoretical optimism}. In other words, we will try to vindicate Gödel's views on phenomenology as a method for arriving at new axioms of ZFC in order to decide independent statements such as CH. Since we have previously written on this very same subject \cite{yoch, yogp1, yogp2}, it is necessary to provide a justification for addressing it once again. 
\\

In \cite{yoch} one can read our first thoughts on this topic. They arose in part as a response against naturalism and ontological objectivism (or realism) in set theory and, unsurprisingly, we tried to defend a pessimistic view, not very far away from that of Hamkins \cite{jdh}. We now find some of its arguments as interesting but badly presented (the dialogue format was certainly not very adequate at that stage of knowledge) and want to partly cover them here. The conclusion at which we arrived was that the overall debate regarding ontological realism and anti-realism was impregnated with a theological flavor. We further developed this thesis in \cite{yodfl}, where we treat Gödel's mathematical intuition, understood as a perception of abstract, spatio-temporal entities, as an example of the mysticism that we wished to attack under the general title of `object-talk'. The conclusion was that the usual mathematical discourse needed to be discharged from all ontological (or psychological) weight, and the adjective that we proposed to describe such a shady region was `liturgical'. 
\\

In \cite{yogp1}, \textit{Gödel's phenomenological program} was first addressed and distinguished from other proposals regarding CH. This also implied a greater recognition of the positive values implicit to Gödel's objectivism. Later, in the second part \cite{yogp2}, the possibility of achieving the evidence of a \textit{dream solution} of CH was recognized from some rather unorthodox views. We now acknowledge a lack of a systematic approach in these papers and how the important place of objectivism in the whole Gödelian project was misunderstood; the main point at which we arrived was that phenomenology was a good candidate for providing a more precise account of what was intended to fall under the loose terminology from \cite{yodfl}.
\\

Here, our desire is to distinguish clearly between competing positions regarding the status of CH and to present Gödel's phenomenological program by introducing a special form of epistemological realism that, in our own opinion, is compatible with both the Gödelian insights and the main methodological assumptions belonging to a general view of phenomenology\footnote{I wish to thank R. Asensi and G. Cobo for their comments and suggestions.}. 
\\

\newpage
\tableofcontents
\newpage

\section*{I. Gödel's program}

In this section we present the usual primary bibliographical sources regarding Gödel's program.

\paragraph{1. Original formulation.} A good place to start at is the \textit{locus classicus} of Gödel's program, namely, his Cantor paper \cite{godelch1}. There, as it is widely known, Gödel defends the view that the Continuum Hypothesis is a proper mathematical problem and thus susceptible of solution. He does so in two steps. First, he motivates a form of ontological objectivism and, from such standpoint, he claims as legitimate the search for new axioms of set theory. The idea then is to establish some criteria in order to accept candidates as legitimate axioms, i.e. as \textit{natural} ones. 

\subparagraph{a. Objectivism.} As early as in his Russell paper \cite{godelrussell}, Gödel defends a strong objectivism regarding sets and concepts. Since the distinction between both notions (see, e.g. \cite{wang}) lies outside the scope of the present paper, let us simply single out the relevant consequences of such view. In few words, the point that Gödel wishes to make is that mathematics bear a \textit{real content} which cannot be `eluded' in any case (say, by a nominalistic conception of sets, see \cite{godelrussell}) and that this content depends upon the conception of sets as objectively existent entities -- even if it is also a possibility to take such existence as a mere working hypothesis (we omit the comparisons with physical objects; see, e.g., \cite{maddy}). Therefore, sets admit a concrete interpretation as `pluralities of things' or as `structures consisting of a plurality of things' (see \cite{godelrussell}). This train of thought is followed later in the Cantor paper:
\begin{quote}
    1.a.1 [...] on the basis of the point of view here adopted, a proof of the undecidability of Cantor’s conjecture from the accepted axioms of set theory (in contradistinction, e.g., to the proof of the transcendency of $\pi$) would by no means solve the problem. For if the meanings of the primitive terms of set theory are accepted as sound, it follows that the set-theoretical concepts and theorems describe some well-determined reality, in which Cantor’s conjecture must be either true or false. Hence its undecidability from the axioms being assumed today can only mean that these axioms do not contain a complete description of that reality. \cite{godelch1}
\end{quote}

The aforementioned point of view consists in an iterative conception of sets, in which one starts with some already well defined \textit{Urelementen} (e.g. the integers) and then iterates the operation `set of'. This is connected to the previous conception of sets \textit{qua} pluralities in the following footnote:

\begin{quote}
    1.a.2 The operation ``set of x's'' cannot be defined satisfactorily (at least in the present state of knowledge), but only be paraphrased by other expressions involving again the concept of set, such as: ``multitude of x's'', ``combination of any number of x's'', ``part of the totality of x's''; but as opposed to the concept of set in general (if considered as primitive) we have a clear notion of this operation. \cite{godelch1}
\end{quote}

Again, we will not discuss here the precise relationship between `set of' and `the concept of set'. Let us simply consider the iterative conception as given. That the present axioms of ZFC do not contain the desired `complete description' of the set theoretical realm mentioned above is a consequence of the following: 

\begin{quote}
    1.a.3 [...] the axioms of set theory by no means form a system closed in itself, but, quite on the contrary,
the very concept of set on which they are based suggests their extension by new axioms which assert the existence of still further iterations of the operation “set of”. \cite{godelch1}
\end{quote}

That is, Gödel's strong objectivism leads to the possibility of solution of CH, that is, to recognize its status as a proper problem\footnote{The precise nature of such problem is a matter of further discussion: see below, section III.}. 

\subparagraph{b. Naturalness.} Thus, we are led to Gödel's idea of \textit{naturalness}. It is worth noting that he had already introduced this very same notion in \cite{godel1938}, even if he would later reject the adoption of the axiom sanctioned in this particular case:

\begin{quote}
    1.b.1 The proposition $A$ [the constructibility axiom $V = L$] added as a new axiom seems to give a \textit{natural} [emphasis is mine] completion of the axioms of set theory, in so far as it determines the vague notion of an arbitrary infinite set in a definite way. \cite{godel1938}
\end{quote} 

In the Cantor paper, however, Gödel refines this idea and splits it in two halves. On the one hand, he first introduces \textit{intrinsic} naturalness:

\begin{quote}
    1.b.2 [...] not only that the axiomatic system of set theory as known today is incomplete, but also that it can be supplemented without arbitrariness by new axioms which are only the natural continuation of the series of those set up so far. [...] there may exist [...]  other (hitherto unknown) axioms of set theory which a more profound understanding of the concepts underlying logic and mathematics would enable us to recognize as implied by these concepts. \cite{godelch1}
\end{quote}

In other words, some axioms may be \textit{semantically implied} by `the concepts underlying logic and mathematics' and hence they will bear some degree of intrinsic naturalness or necessity. This possibility would still be recognized by Gödel in the later \cite{godel1972}: 

\begin{quote}
    1.b.3 [...] there \textit{do} exist unexplored series of axioms which are \textit{analytic} [emphasis is mine] in the sense that they only explicate the content of the concepts occurring in them, e.g., the axioms of infinity in set theory, which assert the existence of sets of greater and greater cardinality or of higher and higher transfinite types and which only explicate the content of the general concept of set. These principles show that ever more (and ever more complicated) axioms appear during the development of mathematics. For, in order only to understand the axioms of infinity, one must first have developed set theory to a considerable extent. \cite{godel1972}
\end{quote}

On the other hand, Gödel presents \textit{extrinsic} naturalness in the Cantor paper as follows:

\begin{quote}
    1.b.4 [...] even disregarding the intrinsic necessity of some new axiom, and even in case it has no intrinsic necessity at all, a probable decision about its truth is possible also in another way, namely, inductively by studying its [...] fruitfulness in consequences, in particular in “verifiable” consequences, i.e., consequences demonstrable without the new axiom, whose proofs with the help of the new axiom, however, are considerably simpler and easier to discover, and make it possible to contract into one proof many different proofs. [...] There might exist axioms so abundant in their verifiable consequences, shedding so much light upon a whole field, and yielding such powerful methods for solving problems [...] that, no matter whether or not they are intrinsically necessary, they would have to be accepted at least in the same sense as any well-established physical theory. \cite{godelch1}
\end{quote}

In this case, axioms are sanctioned as natural on pragmatic reasons, regarding practical fruitfulness. Even if we will restrict ourselves to the analysis of intrinsic criteria, we should bear in mind that the distinction between both kinds of criteria seems both rigid and fundamental to Gödel. The overall program then consists in finding new, either intrinsically or extrinsically, natural axioms of set theory that, moreover, allow us to decide CH. 
\\

\paragraph{2. Philosophical refinements.} It is our belief that Gödel can be read as refining the philosophical tenets of his program in at least two different places: the second version of his Cantor paper \cite{godelch2} and some unpublished material for a lecture \cite{godel1961}. In the former, Gödel wishes to provide a more convincing form of his mathematical objectivism, while in the latter he discusses the general phenomenon of undecidability. 

\subparagraph{a. Refining naturalness.} Let us begin with \cite{godel1961}. Here, Gödel discusses the phenomenon of undecidability in general. For instance, against Hilbert's failed attempt to secure the certainty of mathematics, Gödel argues that 
\begin{quote}
    2.a.1 [...] the certainty of mathematics is to be secured not by proving certain properties by a projection onto material systems - namely, the manipulation of physical symbols but rather by cultivating (deepening) knowledge of the abstract concepts themselves which lead to the setting up of these mechanical systems, and further by seeking, according to the same procedures, to gain insights into the solvability, and the actual methods for the solution, of all meaningful mathematical problems. \cite{godel1961}
\end{quote}

Of course, CH is one of these problems. Hilbert's original intentions (e.g. his axiom of solvability that there is no \textit{ignorabimus} in mathematics) receive then a renovated interpretation in Gödel's thought. This fragment is not only connected to the knowledge required for arriving at intrinsically natural axioms (1.b.2) but also with an open-ended task briefly described in the Russell paper \cite{godelrussell}. There, against Russell's no-class approach, Gödel argues that:

\begin{quote}
    2.a.2 This seems to be an indication that one should take a more conservative course, such as
would consist in trying to make the meaning of the terms “class” and “concept” clearer,
and to set up a consistent theory of classes and concepts as objectively existing entities.
[...] Many symptoms show only too clearly, however, that the primitive concepts need further elucidation. \cite{godelrussell}
\end{quote}

Hence, we see that this pursue of an analysis of our primitive or fundamental notions is ubiquitous in Gödel's thought. This fact concedes even greater importance to answering the following question:
\begin{quote}
    2.a.3 In what manner, however, is it possible to extend our knowledge of these abstract concepts, i.e.,
to make these concepts themselves precise and to gain comprehensive and secure insight into the fundamental
relations that subsist among them, i.e., into the axioms that hold for them? \cite{godel1961}
\end{quote}

Alternatively, this question asks for a way of achieving the desired clarity of the previously mentioned basic notions in order to not only gain understanding of the present axioms but also of new ones which, in turn, will coincide with those that are sanctioned as being intrinsically natural. That is, the notion of intrinsic naturalness would be also further deepened in such case\footnote{With this, we want to express the following. The intrinsic naturalness of some candidate depends on whether it is implied by the meaning of the primitive ideas appearing in it. If we extend our knowledge on these primitive ideas, this would also broaden the derivations that we could made from their meaning and, hence, the domain of intrinsically justified candidates.}. Gödel first notes the following:
\begin{quote}
    2.a.4 Obviously not, or in any case not exclusively, by trying to give explicit definitions for concepts and proofs for axioms, since for that one obviously needs other undefinable abstract concepts and axioms holding for them. Otherwise one would have nothing from which one could define or prove. The procedure must thus consist, at least to a large extent, in a clarification of meaning that does not consist in giving definitions. \cite{godel1961}
\end{quote}

Since we are dealing with the basic concepts, the primitive ideas of the corresponding formal systems, we are not able to establish more axioms or to state some auxiliary results. Thus, there must be a total separation between the purely epistemological and the mathematical tasks that, despite this, must work together in order to decide the initially independent statements formulated in the corresponding formal system: the first provides new axioms and the latter formally settles the question making use of them. Gödel then identifies the epistemological task with phenomenology: 
\begin{quote}
    2.a.5 Now in fact, there exists today the beginning of a science which claims to possess a systematic
method for such a clarification of meaning, and that is the phenomenology founded by Husserl. Here
clarification of meaning consists in focusing more sharply on the concepts concerned by directing our
attention in a certain way, namely, onto our own acts in the use of these concepts, onto our powers in
carrying out our acts, etc. But one must keep clearly in mind that this phenomenology is not a science
in the same sense as the other sciences. Rather it is or in any case should be a procedure or technique
that should produce in us a new state of consciousness in which we describe in detail the basic concepts
we use in our thought, or grasp other basic concepts hitherto unknown to us. I believe there is no
reason at all to reject such a procedure at the outset as hopeless. Empiricists, of course, have the least
reason of all to do so, for that would mean that their empiricism is, in truth, an apriorism with its
sign reversed. \cite{godel1961}
\end{quote}

This is the core of Gödel's conception of phenomenology, as applied to the clarification of intrinsic naturalness\footnote{The reader may point out that, apart the indirect quote from Gödel, we have not yet given any definition of the term `phenomenology'. Moreover, Gödel's general description seems insufficient for providing a robust definition. For the conception of phenomenology that we will defend here, we redirect the skeptic reader to \S5. A preliminary -- though inexact -- explanation of what we understand by phenomenology is, perhaps, this: 
a radical form of conceptual analysis, not far away from the Wittgensteinian analysis of language.}. 

\subparagraph{b. Refining objectivism.} Let us turn our attention to \cite{godelch2}. Here, Gödel vindicates once again his ontological objectivism, now as a response to a comparison with Euclid's fifth postulate and against the claim that the problem for the truth of CH disappears with its undecidability (see \cite{godelch2}). He provides two different kinds of reasons, mathematical and epistemological. On the mathematical side, Gödel argues, there is a far greater asymmetry between having CH or having its negation and the case of geometry and the parallel postulate. In fact, 
\begin{quote}
    2.b.1 The same asymmetry also occurs on the lowest levels of set theory, where the
consistency of the axioms in question is less subject to being doubted by skeptics. \cite{godelch2}
\end{quote}

Regarding the epistemological question,
\begin{quote}
    2.b.2 [...] by a proof of undecidability a question loses its meaning only if the system of axioms under consideration is interpreted as a hypothetico-deductive system, i.e., if the meanings of the primitive terms are left undetermined. \cite{godelch2}
\end{quote}

It is now when Gödel tries to provide new grounds for his strong objectivism. He wants to show how we can, in fact, provide an interpretation of the primitive terms of set theory. His proposal in the notion of mathematical intuition:

\begin{quote}
    2.b.3 [...] we do have something like a perception also of the objects of set theory, as is seen from the fact that the axioms force themselves upon us as being true. I don't see any reason why we should have less confidence in this kind of perception, i.e.. in mathematical intuition, than in sense perception, which induces us to build up physical theories and to expect that future sense perceptions will agree with them, and, moreover, to believe that a question not decidable now has meaning and may be decided in the future. \cite{godelch2}
\end{quote}

Therefore, 1.b.2 may be understood as providing `new mathematical intuitions leading to a decision' of certain statements of set theory such as CH (see \cite{godelch2}) and intrinsic naturalness may be simply described as this `forcing themselves upon us' of certain axioms. Gödel then dedicates a whole paragraph to a further exposition of mathematical intuition (with a terminology heavily inspired by Husserl's). All these efforts, remember, are directed towards the justification of ontological objectivism. Nevertheless, Gödel suddenly introduces a form of \textit{epistemological realism}\footnote{We are following the terminology from \cite{hauser}. Epistemological objectivism or realism differs from the ontological one in that the former gives truth the central place and the latter does so with objects. Wang \cite{wang} recognizes that Gödel would perhaps be also sympathetic towards this, apparently more moderate, form of realism or objectivism in mathematics.}:

\begin{quote}
    2.b.4 However, the question of the objective existence of the objects of mathematical intuition [...] is not decisive for the problem under discussion here. The mere \textit{psychological fact} [emphasis is mine] of the existence of an intuition which is sufficiently clear to produce the axioms of set theory and an open series of extensions of them suffices to give meaning to the question of the truth or falsity of propositions like Cantor’s continuum hypothesis. What, however, perhaps more than anything else, justifies the acceptance of this criterion of truth in set theory is the fact that continued appeals to mathematical intuition are necessary [...] for obtaining unambiguous answers to the questions of transfinite set theory [...]. This follows from the fact that for every axiomatic system there are infinitely many undecidable propositions of this type. \cite{godelch2}
\end{quote}

This is, perhaps, one of the most intriguing remarks ever made by Gödel. He leaves us with the following question: which could be the nature of this `psychological fact'?

\newpage

\section*{II. Gödel's global phenomenological program}

We believe that a special form of epistemological objectivism is suitable for reconciling Gödel's views with Husserl's original demands. This belief is derived from: (i) a vindication of phenomenology in the overall program as conceived by Gödel and (ii) a complete separation between the phenomenological and the mathematical in such program. In this section we will make use of secondary bibliographical sources, notoriously, \cite{wang}.

\paragraph{3. The phenomenological.} We believe that the adjective `psychological' in 2.b.4 should be better understood as \textit{phenomenological}. Indeed, it is highly unlikely that Gödel is using this adjective in order to describe something close to a subjective experience. Rather, the power of mathematical intuition as some kind of `perceiving' applicable to abstract notions is precisely that it is independent from the psychical. 

\subparagraph{a. Against self-evidence.} As a response to our view, one can perhaps invoke the following fragment from \cite{godel1961}:
\begin{quote}
   3.a.1 [...] one has examples where, even without the application of a systematic and conscious procedure [i.e. the phenomenological method], but entirely by itself, [...] it turns out that in the systematic establishment of the axioms of mathematics, new axioms, which do not follow by formal logic from those previously established, again and again become evident. It is not at all excluded by the negative results mentioned earlier that nevertheless every clearly posed mathematical yes-or-no question is solvable in this way. For it is just this becoming evident of more and more new axioms on the basis of the meaning of the primitive notions that a machine cannot imitate. \cite{godel1961}
\end{quote}
It is natural to compare this to 2.b.4 by establishing an analogy between the `psychological fact' and the axioms that become evident `by themselves'. Under this reading, Gödel is then saying that the phenomenological method, so dramatically presented in 2.a.5, would only be one possibility in searching intrinsic justifications for new axioms and that there are, still, \textit{pure psychological} phenomena of evidence, irreducible to systematic analyses. 
\\

We interpret this, rather, as Gödel simply trying to make the problem of intrinsic justification attractive by itself, i.e. without appealing to the phenomenological method. Indeed, in \cite{godel1961}, he later returns to phenomenology in connection with Kant's philosophy of mathematics. Moreover, leaving aside bibliographical questions, it is clear to us that we do not have any way of studying how axioms are made evident `by themselves' or `force themselves upon us' (2.b.3). If we, without further efforts, embrace this shaky ground of self-evidence of axioms and identify it with intrinsic justification, nothing prevent us from deeming the intrinsic approach in Gödel's program as limiting or even purely subjective \cite{bagaria, fefermannewaxiomsmaddy, fontanella}.

\subparagraph{b. Wang's examples.} In \cite[8.2]{wang}, intrinsic justifications are provided for some already accepted axioms of set theory, AI (\textit{infinity}), AE (\textit{extensionalization}), AF (\textit{foundation}), AS (\textit{separation}), AP (\textit{powerset}) and AR (\textit{replacement}), in connection to various Gödel's remarks. Let us briefly go through them in order to analyze whether the presence of the phenomenological method among the intrinsic justifications is merely accidental. The basic tools required for such justifications are \textit{overview}, \textit{idealization} and \textit{extensionalization} (see \cite{wang}). The adjective `psychological' appears crucially in the first of the following fragments (reporting two observations made by Gödel) illustrating how we can overview sets and, through idealization, form the subset of any given set:
\begin{quote}
    3.b.1 To arrive at the totality of integers involves a jump. Overviewing it presupposes
an [Wang: idealized] infinite intuition. In the second jump we consider not only
the integers as given but also the process of selecting integers as given in intuition.
``Given in intuition'' here means [Wang: an idealization of] concrete intuition. Each
selection gives a subset as an object. Taking all possible ways of leaving elements
out [Wang: of the totality of integers] may be thought of as a \textit{method} for producing these
objects. What is given is a \textit{psychological analysis} [emphasis is mine], the point is whether it produces
objective conviction. This is the beginning of analysis [Wang: of the concept of set]. \cite[7.1.18, cf. 8.2.11]{wang}
\end{quote}

\begin{quote}
    3.b.2 We idealize the integers (a) to the possibility of an infinite totality, and (b) with omissions. In this way we get a new concretely intuitive idea, and then one goes on. [...] What this idealization --realization of a possibility-- means is that we conceive and realize the possibility of a mind which can do it. We recognize possibilities in our minds in the same way as we see objects with our senses. \cite[7.1.19]{wang}
\end{quote}
According to Wang, these remarks intrinsically justify both AI, for we can recognize that there are infinite sets by the very nature of mathematical intuition, and AP, since we have a `method' for forming subsets: 
\begin{quote}
    3.b.3 The axiom of subset formation comes before the axiom of power set. We can form the power set of a set, because we understand the selection process [the `method' above] (of singling out any subset from the given set) intuitively, not blindly. \cite[8.2.17]{wang}
\end{quote}
According to Wang, however, extensionalization is also relevant in this case:
\begin{quote}
    3.b.4 The step from the existence of all subsets of an infinite set to the overview ability of its power set clearly involves a strong idealization of our intuition. This matter of presupposition, so far as existence is concerned, is not a question of temporal priority. The point is, rather, that, conceptually, objects have to exist in order for the set of them -- as their unity-- to exist. \cite{wang}
\end{quote}
In turn, AE holds because it is simply a defining property of sets, i.e. because one basic feature of sets is extensionalization, while AF is justified by the non-artificiality of the rank hierarchy (see \cite[8.7]{wang}). Wang also provides a simple justification of AR by using overview. Notoriously, Gödel seemed to differ with it since, according to him and differently from the other axioms, AP lacks `immediate evidence', i.e. `previous to any closer analysis of the iterative concept of set' (see \cite[8.2.14]{wang}).
\\

Now, we believe that all of these justifications can be seen as legitimate and valid examples of applications of the phenomenological method. What both uses of the adjective `psychological', in 2.b.4 and 3.b.1, have in common is that they are related to some properties of mathematical intuition, be it regarding its existence in certain contexts, as in the first case, or having to do with its idealization, as in the latter. But mathematical intuition is not a psychological faculty. The phenomenological reduction suspends our judgment regarding positive knowledge, psychological considerations included. When Husserl talks about extending the way in which we usually talk about perception, including grasping essences, he does so in a non-psychological manner (see, e.g., \cite{scherer}). The kind of relation that Husserl wishes to work out is \textit{ideal}, as confronted to the \textit{real}. It is then consistent, from what Gödel explains in \cite{godelch2} and his conception of phenomenology (see \cite{wang}), to regard mathematical intuition as providing ideal contents, not real. And psychological contents certainly fall under the latter kind\footnote{Despite all of this one could, after all, reject any need for making use of the notion of intuition. A very reasonable alternative consists in simply leaving these psychological or `internal' considerations aside and placing our attention on the practical, historical or `external' factors. Thus, it can be argued that the usual extensions of ZFC gain or will gain weight by considerations of the latter kind and that every other purported justification will actually be an \textit{a posteriori} ornament. The problem, however, is that accounts of this form usually make use of terms like `foundation' or `justification'. We present this problem and its subsequent reply in \S7.d and \S8.}. 

\subparagraph{c. Concrete criteria.} The previous intrinsic justifications were given for already known axioms of set theory. Some doubts may arise regarding a similar justification for, say, AC. However, the mere necessity for the search of new axioms that is provided by Gödel's program is enough to question how we can extend intrinsic justifications, not for dealing with \textit{actual} statements but for \textit{potential} ones. Wang distinguishes here two different problems: 
\begin{quote}
    3.c.1 What are the principles by which we introduce the axioms of set theory? This is different from the related question, What is the precise meaning of the principles, and why do we accept them? [Gödel] did not say much about the second question, but he seemed to suggest that it should be answered by phenomenological investigations in the manner of Husserl. \cite[8.7]{wang}
\end{quote}
Taking this at face value would, as in the previous case of self-evidence, compromise the relevance of phenomenology for the whole of Gödel's program. Indeed, if we could simply recognize several principles regulating the introduction of new natural axioms then, once again, the method described in 2.a.5 as an answer to 2.a.3 would be superficial. If we do not have \textit{grounds} for the adoption of such principles, then it is safe for us to simply take them as mere contingencies of our practice. 
\\

For a clear exposition on the development of Gödel's views regarding the regulating intrinsic criteria, as presented in \cite[8.7]{wang}, we redirect the reader to \cite{bagaria}. We shall make, however, some observations. First, it is convenient to note that they are not mutually exclusive, i.e. the same axiom may be sanctioned as natural by two different criteria (see \cite[8.7.6]{wang}). On the other hand, at some point Gödel believed that 
\begin{quote}
    3.c.2 All the principles for setting up the axioms of set theory should be reducible to a form of Ackermann's principle: The Absolute is unknowable. The strength of this principle increases as we get stronger and stronger systems of set theory. The other principles are only heuristic principles. Hence, the central principle is the reflection principle, which presumably will be understood better as our experience increases. Meanwhile, it helps to separate out more specific principles which either give some additional information or are not yet seen clearly to be derivable from the reflection principle as we understand it now. \cite[8.7.9]{wang}
\end{quote}

What we want to bring to our attention here is the following. Wang argues that, from Ackermann's principle, Gödel is able to \textit{derive} three special principles, which are to be considered as a considerable refinement of his previous attempts in providing a comprehensive list of intrinsic criteria (see \cite[8.7.10]{wang}). Now, a natural reading of the preceding fragment is that, despite Wang's initial comment 3.c.1, Gödel did actually make some progress in the direction of the second question, namely, he claimed that Ackermann's principle may help to justify the corresponding regulating principles or, in other words, that it may help us as a definite \textit{grounding} for the acceptance of such principles. It is clear that such principle is intended as a further elaboration on the analyses of the primitive concepts of set theory -- following 1.b.2 and 2.a.3 -- and that, therefore, it provides additional insights on the behavior of such meanings.

\subparagraph{d. Further comments.} Since we have been comparing the subjective with the nature of mathematical intuition, it seems reasonable to take a look at some (rather elusive, Wang confesses) comments made by Gödel on the first. Wang relates these observations with the Cantor-von Neumann's axiom: a collection of sets is a set if and only if is not as large as $V$ (see \cite[8.3]{wang}). Following 3.b.1 and 3.b.2, Gödel establishes some limitations regarding the overviewing of sets: 

\begin{quote}
    3.d.1 This significant property of certain multitudes [that they can be overviewed] must come from some more solid foundation than the apparently trivial and arbitrary phenomenon that we can overview the objects in each of these multitudes. [..] without the objective picture, nothing seems to prevent us from believing that every multitude can be thought together. [...] Some pluralities can be thought together as unities, some cannot [due to the contradictions]. Hence, there must be something objective in the forming of unities. \cite[8.3.1]{wang}
\end{quote}

Or, more precisely, he distinguishes the mere psychological act and a strong objective sense of overviewing sets. The difference, as Wang argues, is that between mere (subjective) \textit{thought} and (objective) \textit{intuition}, and this induces a subjective and an objective picture, or view, of sets, respectively: 

\begin{quote}
    3.d.2 In some sense, the subjective view leads to the objective view. Subjectively, a set is something which we can overview in one thought. [...] Different persons can, we believe, each view the same multitude in one thought. Hence, it is natural to assume a common nucleus which is the objective unity. [...] Idealization is decisive in both cases [Wang: the subjective and the objective unity]. \cite[8.3.2]{wang}
\end{quote}

The last line of this paragraph can be contextualized with the help of the following, which contains another notorious occurrence of the adjective `psychological': 

\begin{quote}
    3.d.3 [...] the range of possible knowledge is wider than the range of existence that can be justified from the subjective viewpoint. The \textit{psychological act} [emphasis is mine] of thinking together all objects of a multitude in one thought yields more sets from the objective viewpoint because stronger idealizations are appropriate [Wang: than from the subjective viewpoint]. From the idealized subjective view, we can get the power set. But the indefinability of the universe $V$ of all sets can't be got by the subjective view at all. The difference in strength becomes clear only when you introduce new principles which make no sense at all in the subjective view. For every set there is some mind which can overview it in the strictest sense. \cite[8.3.3]{wang}
\end{quote}

Roughly, as Wang concludes, thought can go beyond intuition; idealization can extend intuition with the help of thought (see \cite[8.3]{wang}). There is a crucial detail in this fragment, namely, that Gödel recognizes a psychological act `from the objective viewpoint'. This may be interpreted in many ways, but we believe that it is consistent with Husserl's overall project of phenomenology to clarify it as implying a clear-cut distinction between the merely subjective and the \textit{systematic subjective}. Even so, the `psychological' must be dropped in this context and, instead, we argue again, it may be rephrased as `phenomenological'. Additionally, Gödel seems to distinguish `idealized' subjective and objective views, perhaps meaning that the application of idealization can be carried out in two different ways or, better, differently in both the subjective and the objective pictures\footnote{It is also likely that Gödel is having here in mind the noetic-noematic distinction, which is essential for phenomenology. Properly explaining this sketch of an idea requires, of course, further work.}. In any case, the former does not seem apt in order to explain the nature of intrinsic criteria, since there are principles that `make no sense at all' from such viewpoint, Ackermann's principle in particular. From our perspective, the adjective `psychological' applies to the subjective and idealized subjective mindsets, while its (improper) application for the case of the objective view should be better understood as describing the domain of the phenomenological\footnote{There are additional sources for our reading of Gödel. For example, one could connect our comments above with Gödel's thought on computabilism and the existence of some systematic procedures that fail to be implemented in a Turing machine. Additional fragments can be read in \cite{wang} regarding a formal system and its extensions and how mathematical intuition is essentially required to do so; a computer is not able to achieve such tasks (cf. 3.a.1). This is consistent with our treatment in \cite{yowitt} of Kripke's skeptic paradox regarding rule-following. There, we argued that, given a definite amount of pluralism applied to the potential extensions of a formal system, there is no rational justification (apart from tradition and historical force) for adopting one against the other. Of course, mathematical intuition as understood by Gödel undermines pluralism and thus makes possible adopting a concrete extension of a formal system as following \textit{naturally} from our original intentions given in the initial formalization of the primitive notions.}. 
\\

\paragraph{4. The phenomenological and the mathematical.} The aim of the preceding paragraph was to vindicate the role of phenomenology in Gödel's program. To put it in a few lines: phenomenology --  against self-evidence and  psychological explanation -- should be the preferred tool for clarifying the primitive notions of set theory and, thus, for providing insights and analyses in order to elaborate intrinsic justifications of new axioms. The point that we want to make now, however short, is that this phenomenological side of Gödel's program or, as we prefer to name it, \textit{Gödel's phenomenological program}, should be distinguished from the purely mathematical one, namely, the actual derivation of the desired statements in the natural extension of ZFC provided by the former. 
\\

We claim that, besides 2.a.4, one can read Gödel as separating both tasks in other fragments reported in \cite{wang}. For example, the following fragment seems to confirm our view, at least understanding `epistemological' to mean `phenomenological' in this context, as we wish that the preceding discussion has helped to convey (see our comments above between 2.a.4 and 2.a.5):
\begin{quote}
    4.a.1 The epistemological problem is to set the primitive concepts of our thinking right. For example, even if the concept of set becomes clear, even after satisfactory axioms of infinity are found, there would remain more technical (i.e., mathematical) questions of deciding the continuum hypothesis from the axioms. This is because epistemology and science (in particular, mathematics) are far apart at present. It need not necessarily remain so. True science in the Leibnizian sense would overcome this apartness. In other words, there may be another way of analyzing concepts (e.g., like Hegel's) so that true analysis will lead to the solution of the problem. \cite[7.3.21]{wang}
\end{quote}

\begin{quote}
    4.a.2 This scarcity of results, even as to the most fundamental questions in this field [such as CH], may be due to some extent to purely mathematical difficulties; it seems, however (see Section 4 below), that there are also deeper reasons behind it and that a complete solution of these problems can be obtained only by a more profound analysis (than mathematics is accustomed to give) of the meanings of the terms occurring in them (such as ``set'', ``one-to-one correspondence'', etc.) and of the axioms underlying their use. \cite{godelch1}
\end{quote}

Even if the difference between the epistemological and the mathematical is due to the insufficiency of the present knowledge, such duality persists in our case and is, for now, a problem that is to be left aside. Another negative consideration is the next one: 
\begin{quote}
    4.a.3 At present we possess only subjective analyses of concepts. The fact that such analyses do not yield decisions of scientific problems is a proof against the subjectivist view of concepts and mathematics. \cite[7.3.22]{wang}
\end{quote}
This, together with Gödel's general attitude towards phenomenology \cite{hauser, tieszen, wang} seems to indicate that a desideratum of the phenomenological method, as applied to such clarification of concepts, would consist in transcending mere subjective analyses or, alternatively, that it would emerge as a definite tool in order to go beyond these and pursue objective ones. 
\\

\paragraph{5. The phenomenological, revisited.} For the expert phenomenologist, the statement above about how one should read Gödel's expression `psychological fact' as `phenomenological fact' will surely appear intriguing and in need of elucidation. Indeed, we have made several comments regarding the phenomenological, its domain and nature, but we have not settled a basis for the reader to systematically compare our view with others. We believe that the correct way to proceed now is to provide several examples or insights illustrating the nature of phenomenology. We hope that these brief observations will make our exposition more clear and exact.

\subparagraph{a. Anti-reductionism.} A notorious example of a discussion of phenomenology in the context of the philosophy of mathematics is given by Rota \cite{rota}. Against the usual metaphysical debates regarding mathematical objects, Rota believes that other questions, close to the tradition of phenomenology, are of more importance (see, e.g. \cite[VII, XII]{rota}). In \cite[XV]{rota}, Rota argues that the Husserlian notion of \textit{Fundierung} (roughly, \textit{foundation}) should be vindicated as a logical notion. The relevance of this concept for mathematics can be seen in the following excerpt:

\begin{quote}
    5.a.1 [...] the relation between the truth of mathematics and the axiomatic truth which is indispensable in the presentation of mathematics is a relation of \textit{Fundierung}. \cite[XI]{rota}
\end{quote}
In our present context, this remark contributes to the demarcation between the epistemological and the mathematical. Rota motivates the analysis of \textit{Fundierung} by a series of examples regarding the distinction between `function' and `facticity'. Through these terms, Rota establishes the main feature of anti-reductionism in phenomenology:

\begin{quote}
    5.a.2 \textit{Fundierung} is a primitive relation, one that can in no way be reduced to simpler (let alone to any ``material'') relations. It is the primitive logical notion that has to be admitted and understood before any experimental work on perception is undertaken. Confusing function with facticity in a \textit{Fundierung} relation is a case of \textit{reduction}. Reduction is the most common and devastating error of reasoning in our time. Facticity is the essential support, but it cannot upstage the function it \textit{founds}. \cite[XV]{rota}
\end{quote}

We believe, with Rota, that \textit{Fundierung} can be explained by establishing some analogies between Husserl's and Wittgenstein's thought (again, see \cite[XV]{rota}), namely, as essentially corresponding to the one existing between an element (`facticity') and its role (`function') within a Wittgensteinian language game (see \cite{wittPI} and \cite[XV]{rota})\footnote{Hence, for the reader that wishes to omit the notion of \textit{Fundierung}, we can simply say that one of the tasks of phenomenology is to find relations of the type facticity-function while properly distinguishing their two poles. This is, in turn, the anti-reductionistic character of phenomenology that we wish to establish here.}. 

\subparagraph{b. Immediacy.} A happy coincidence in the history of philosophy is that Wittgenstein dealt during some time with the idea of a `phenomenological language' (see \cite{wittphenomenology}). Of course, it would be lazy to take Wittgenstein's remarks and attach them, without any critical remorse whatsoever, to Husserl's. Nevertheless, since we are providing here an interpretation of phenomenology, this happy coincidence allows us to let Wittgenstein speak for ourselves.
\\

As Wang recognizes, both Husserl and Wittgenstein seek to reduce philosophical problems to some kind of immediacy or `fundamental intuition', while the difference between them is \textit{where} they intend to locate it (see \cite[10.1]{wang}). Husserl, on the one hand, speaks about intentionality, about noema and noesis, about reflection of our various \textit{Erlbenisse} and consciences. Rota also follows Husserl when speaking about evidence as a logical notion (see \cite[XIV]{rota}). Wittgenstein, on the other hand, `begins and ends with the perceptual immediacy of our intuition of the actual use of words in a given situation' (\cite[10.1]{wang}). Indeed, we can read:

\begin{quote}
    5.b.1 The investigation of the rules of the use of our language, the recognition of these rules, and their clearly surveyable representation amounts to, i.e. accomplishes the same thing as, what one often wants to achieve in constructing a phenomenological language.
    
    Each time we recognize that such and such a mode of representation can be replaced by another one, we take a step toward that goal.
    
    [...] What we’re missing isn’t a more precise scrutiny [...] nor the discovery of a process behind the one that is observed superficially (that would be the investigation of a physical or psychological phenomenon), but clarity in the grammar of the description of the old phenomenon. Because if we looked more closely we would simply see something else, and would have made no advance on our problem. \textit{This} experience, and not another, is what needs to be described. \cite{wittphenomenology}
\end{quote}

From our point view, this establishes a strong analogy between the idea of the \textit{grammatical} in Wittgenstein and the \textit{ideal} in Husserl's Prolegomena \cite{husserlILI}. Rota's essential notions of \textit{Fundierung} and evidence are ideal, or grammatical, because they cannot be reduced to the psychological, physical or metaphysical, but they share a certain degree of immediacy. Psychological explanations are concerned with the \textit{real} and phenomenology is occupied with the \textit{ideal}, only that not in an idealist or metaphysical sense, but in a radically grammatical one. Of course, immediacy is not to be confused with ecstatic receptivity of data as given. Rather, we should profess an active search for new `modes of presentation'. Following Rota we may say that, where Wittgenstein proposes, for such purpose, an analysis of language by detecting language games and the use of certain terms, Husserl favors phenomenological reduction and eidetic variations with the eye on essences (see \cite[XV]{rota})\footnote{Therefore, we may further elaborate on our preferred notion of phenomenology (see the footnote above) by establishing it as a kind of language analysis. The immediacy feature is given here by our ability to recognize the language-game relations described before. This puts more emphasis on how the ideal/grammatical mentioned in the text stands far away from the psychological.}.

\subparagraph{c. Descriptivity.} As Scheler believed, phenomenology is best described as an attitude towards philosophical problems, instead of a science with highly restrictive methodological standards or a rigid method. One strong reason in doing so is that phenomenology has not experienced a similar development to other sciences, nor it seems to be able to achieve it in the future. Nevertheless, we may still describe it as a \textit{scientific attitude towards} these problems, i.e. with certain degree of rigor\footnote{One may argue that `having faith in' is also to count as a possible attitude towards philosophical problems. What we argue here is that phenomenology differs with it in its scientific character (against a supposedly dogmatic one).}. As Rota recognizes, we cannot hope to provide a clear-cut definition of, say, \textit{Fundierung} or a mathematical formalization of it:
\begin{quote}
    5.c.1 The impossibility of formalization should not be confused with lack of rigor. Formal presentation is not the only kind of rigor. In philosophy, presentation by examples is an essential element of rigor. Examples are to philosophical discourse what logical inference is to mathematical proof. \cite[XV]{rota}
\end{quote}

This is similar in spirit to \cite{husserlFCE}, where Husserl presents phenomenology as an strict science, even though its criteria of rigor differ from those of natural sciences. If we follow Wittgenstein (and the early Husserl), it seems that phenomenology is to be understood as a systematic process, which is nevertheless descriptive, i.e. neither positive or negative:

\begin{quote}
    5.c.2 [...] it actually isn’t a question of the difficulty of calling up
a particular image before my mind’s eye, nor is it a question of something that I can try but fail at; rather,
it’s a question of acknowledging a rule for a mode of expression. \cite{wittphenomenology} 
\end{quote}

Or, perhaps more crudely put:

\begin{quote}
    5.c.3 Philosophy must not interfere in any way with the actual use of language, so it can in the end only describe it.
    
For it cannot justify it either.

It leaves everything as it is.

It also leaves mathematics as it is, and no mathematical discovery can advance it. A ``leading problem of mathematical logic'' is for us a problem of mathematics like any other. \cite[I, \S 124]{wittPI}
\end{quote}

Now, as Wang also points out, the attitude in Wittgenstein and Husserl towards science differs profoundly (see \cite[10.1]{wang}). The former wishes to extirpate misleading forms of expression that arise in philosophical discourse pretending to imitate science; the latter, to turn philosophy into a rigorous science, where fundamental truths may be obtained. But a candid reading of both could bring the two positions together: if Wittgenstein separates philosophical and scientific discourse is because he has a clear image in mind of the demands that philosophy, as an activity, should meet; if Husserl wishes to label philosophy as a \textit{sui generis} science, it is because he believes that philosophical and (natural) scientific discourse should be radically separated. From this, we deem our position --wobbly standing between both approaches -- as provisionally tenable and as bearing some degree of reasonability.

\subparagraph{d. Openness.} The last and probably most important feature that we want to make explicit regarding our conception of phenomenology is that of openness or \textit{globalness}. This, as a general methodological prescription, amounts to radically eluding any kind of mysticism, i.e. unjustified restrictions of the method itself\footnote{The skeptic reader perhaps wonders if admitting supposedly irreducible and fundamental concepts (such as that of \textit{Fundierung}) is not a sign of mysticism itself. This is connected with our desire of finding in phenomenology a scientific attitude of approaching, asking, sharing, revising and concluding problems, not an empty title under which one is licensed to attain truths without further critique. Phenomenology means for us something close to what philosophy meant for Wittgenstein (even in the \textit{Tractatus}): a continuous activity \cite{wittTLP}.}. Descriptivity, understood as implying the absence of suppositions -- which is ubiquitous in Husserl's phenomenology in virtue of the \textit{epoché} -- leads to extending the scope of our analyses to domains of the greatest generality possible\footnote{By `absence' we only mean that the phenomenological method implies that we act \textit{as if} previous knowledge were inexistent. The point is not to abandon our previously held beliefs (as in the Cartesian strategy of radical doubt) but to methodologically act \textit{as if} these contents were not active. Properly distinguishing the method of phenomenological reduction from the Cartesian one will help to understand the openness that we intend our conception of phenomenology to have.}. 
\\

For example, when describing Gödel's (local) phenomenological program, several limitations seem to arise: (a) regarding the iterative conception of set, together with the theological insights provided by Ackermann's principle, as the unique one and (b) excluding the possibility, in principle, of arriving at new primitive ideas. But, in fact, the essential part of the aforementioned program needs not to respond to these prefigurations and, rather, is able to drop them when demanded by its very own task. As a hint for how these restrictions can be overcome, we may think about tools already provided by Husserl that take the form of \textit{local phenomenological analyses}, i.e. (i) local phenomenological reductions bracketing only concrete domains of knowledge and (ii) local eidetic variations regarding concepts with some already established properties. 
\\

Following this path, we may then clarify intrinsic and extrinsic criteria as follows. The adjective `intrinsic' refers to the domain obtained when operating a local phenomenological reduction with respect to certain notion of set (e.g. the iterative conception); `extrinsic' refers to the domain obtained when operating such reduction but with respect to certain knowledge (e.g. certain consequences of the axiom). The search for an intrinsically justified axiom then consists in realizing an eidetic variation of certain notion of set but in the sphere of the intrinsic; the search for an extrinsically justified axiom consists in realizing an eidetic variation of certain notion of set but in the domain of the extrinsic. \textit{Gödel's global phenomenological program}, when restricted subjected to restraints (a) and (b) yields what we have been labeling as \textit{Gödel's phenomenological program} and, similarly, it may result in different phenomenological programs, each defined by some leading methodological restrictions\footnote{This idea, together with some key comparisons between Husserl and Wittgenstein, is motivated in \cite{yogp1}. See also \cite{yoch} for a presentation of Wittgenstein's conception of a phenomenological language in the context of the problem of CH.}.

\paragraph{6. Gödel's optimism.} Now, our conception of phenomenology is probably insufficient for the Gödelian. The question that remains is: how are we to arrive at precise new axioms by a method that is merely descriptive? 

\subparagraph{a. Optimism.} A key feature belonging to Gödel's thought and lacking in the previous exposition is that of \textit{optimism} (see \cite{wang}): indeed, epistemological realism, as a belief in our capacity in rendering every proper mathematical problem as either true or false, implies the possibility of adding new axioms when applied to undecidable statements. In particular, such optimism translates to a phenomenological one, once we consider this same possibility from the intrinsic point of view -- which, as we have argued all along, must take the form of a phenomenological analysis.
\\

Gödel's optimism, according to Wang, takes many forms in his philosophy, being notorious when related to the phenomenological method, scope and results. For instance, Gödel believed it possible to achieve Husserl's foundations for a systematic monadology \cite{husserlMC}, developed later in the style of Leibniz, through the application of the phenomenological method (see, e.g., \cite[5.3.6, 5.3.7, 5.3.33]{wang}). Despite recognizing that finding `a list of the main categories (e.g., causation, substance, action) and their interrelations' would be an alternative to the phenomenological method, such notions `however, are to be arrived at phenomenologically' (see \cite[5.3.7]{wang}). In general, he seemed to believe that 
\begin{quote}
    6.a.1 Philosophy aims at a theory. Phenomenology does not give a theory. In a
theory concepts and axioms must be combined, and the concepts must be precise
ones. \cite[9.3.10]{wang}
\end{quote}

This, at least in principle, is consistent with our previous remarks on the descriptivity of the phenomenological method (as does the use of `analytic' in 1.b.3). However, the `theory' at which philosophy aims is probably, according to Gödel, a `science of metaphysics' (see \cite[5.3.8]{wang}) and axiomatic in nature (see \cite[9.1]{wang}). Therefore, if the phenomenological method seems powerful enough for Gödel for settling a strict metaphysical system, it is not unsurprising that he also to did so regarding the search of new axioms of set theory: 
\begin{quote}
    6.a.2 [...] This science of intuition is not yet precise, and people cannot learn it yet. At present, mathematicians are prejudiced against intuition. Set theory is along the line of correct analysis. \cite[5.3.29]{wang}
\end{quote}

\subparagraph{b. Sudden illumination.} Now, phenomenology is presumably a descriptive science or method, in which we can only hope to attain a better understanding of the analyzed ideas and, as Gödel wishes, of the primitive ideas of set theory. Our question above remains: we feel that, without the establishment of \textit{synthetic} (i.e. positive or negative) results and facts, we are left with nothing, and that, in case we did foresee something, such product of the method would be, at best, subjective.
\\

Now, Gödel's optimism relies upon a more convoluted component, namely, \textit{sudden illumination} or \textit{conversion} (see \cite[9.1]{wang}). In a nutshell, Gödel believes that, in the process of applying the phenomenological method, we can expect a possible change in our \textit{attitude} (recall 2.a.5) so that, for example, certain situations will seem for us different as they were before. Leaving aside the comments made by Gödel on Husserl's radical change of perspective (see \cite[5.3.30, 31, 32]{wang}), it is interesting to reproduce the following ones:

\begin{quote}
    6.b.1 [...] It is possible to attain a state of mind to see the world differently. One fundamental idea is this: true philosophy is [Wang: arrived at by] something like a religious \textit{conversion} [emphasis is mine]. \cite[9.1.14]{wang}
\end{quote}

\begin{quote}
    6.b.2 [...] This is different from doing scientific work; [Wang: it involves] a change of personality. \cite[9.1.15]{wang}
\end{quote}

Hence, according to Gödel, if we pursue phenomenological analysis, things will then appear to us under a new light. Note, however, that phenomenology -- as some definite activity -- only provides the possibility of a common attitude towards certain problems while Gödel is here establishing his \textit{faith} in that this very same attitude will lead us to new, concrete, advances. In other words, epistemological realism must \textit{precede} the application of the phenomenological method in order to make its success plausible. As far as Gödel goes, such faith is inseparable from his epistemological realism and overall optimism, but we shall put emphasis on the demarcation between the method (i.e. phenomenology) and the faith in it that he seemed to profess (\textit{via} his epistemological realism).

\subparagraph{c. Phenomenological optimism.} In fact, the hypothesis of epistemological objectivism is perhaps too vague and strong for this context. We can embrace the weaker belief in that the corresponding phenomenological analyses of the primitive ideas can always be pursued and gradually achieve new properties and relations held true by them. Indeed, if this is the case, then we can expect to obtain new natural axioms and then try to decide CH from them. Notoriously, this position can be linked with Wittgenstein's `complete analyses' from \cite{wittTLP} (whose approach, incidentally, Gödel seemed to prefer to the later \cite{wittPI}). In fact, one can read Gödel as defending this claim in other places, for example, regarding the specific case of the Church-Turing Thesis (CTT): 
\begin{quote}
    6.c.1 If we begin with a vague intuitive concept, how can we find a sharp concept
to correspond to it faithfully? The answer is that the sharp concept is there all
along, only we did not perceive it clearly at first. This is similar to our perception
of an animal first far away and then nearby. We had not perceived the sharp concept
of mechanical procedures before Turing, who brought us to the right perspective.
And then we do perceive clearly the sharp concept. \cite[7.3.1]{wang}
\end{quote}
The comparison with the perception of the physical should not miscarry us -- we believe -- to read Gödel as embracing some metaphysical position: the conception of phenomenology made explicit above warns us that, indeed, we may improperly speak of such perception but that, at the end, such form of discourse is a non-reductionistic account of how we deal with (abstract) notions, i.e. one that is independent from both the psychological and the metaphysical. Our opinion is that this can be argued on the basis of the following excerpt, which is consistent with our own remarks above:
\begin{quote}
    6.c.2 ``Trying to see (i.e. understand) a concept more clearly'' is the correct way
of expressing the phenomenon vaguely described as ``examining what we mean
by a word.'' \cite[7.3.4]{wang}
\end{quote}

Now, even if carrying out a phenomenological analysis does not consist in a passive activity (like, e.g., waiting for some kind of revelation), Gödel acknowledges that, in the process of analyzing a concept, new properties will became evident and, in this sense, they will be \textit{forced upon us}:
\begin{quote}
    6.c.3 If there is nothing sharp to begin with, it is hard to understand how, in
many cases, a vague concept can uniquely determine a sharp one without even the
\textit{slightest} freedom of choice. \cite[7.3.3]{wang}
\end{quote}
This, of course, without any mystical component, remains as an open possibility. It could be that, when analyzing a concept, we may end up with two different ones because we have been able to extend it in two incompatible ways. But even then we may then try to preserve the common features as the underlying notion to both of them, etc. It is interesting to note how Gödel sees the notion of Turing machine as the sharp notion corresponding to algorithm, i.e. in need of no further analysis, while at the same time considers the notion of set as susceptible of more sharpening. This may yield some properties pertaining that status of CH against that of CTT\footnote{This topic is of course appealing by itself, but ut fails outside the scope of our present considerations. One could try to establish and break the analogy at different places between both statements. For instance, both deal with \textit{informal} notions: CH by relying on finding new axioms that are `natural' and the CTT by establishing the equivalence between an informal notion and its formal correlate. CTT, on the other hand, may be considered as bearing some quasi-empirical nature, while this hardly(?) applies to CH.}.
\\

The main issue here is to recognize in Gödel's optimism a form of \textit{phenomenological optimism}, namely, that we can make concepts gradually clearer through analyses of language and, thus, make new essential (analytic) properties evident. The other half of Gödel's program, namely, the mathematical one, remains open. For Gödel, as we have seen before, the epistemological and the mathematical are differentiated \textit{accidentally} (i.e. given the present state of philosophy with respect to science) while for us such distinction is \textit{essential}. With our phenomenological optimism we wish to emphasize that the new relations holding between the primitive ideas of set theory come first, and that the mathematical task of determining whether these settle CH is a secondary task, at least from the intrinsic point of view. In a nutshell, Gödel's epistemological realism entails that CH has a solution and, in particular, that new natural axioms will be made explicit (but with an eye put on deciding CH). Our phenomenological realism implies the latter and leaves the former open, but we recognize through it that a solution is \textit{possible} and, hence, we stand optimistic with respect to the status of CH. As we have previously seen, epistemological realism implies phenomenological optimism, but the converse is not true, at least in principle, since it could be the case that all that could be intrinsically justified as a new axiom would not decide CH, for example.

\newpage 
\section*{III. Discussion}

In this section we confront our phenomenological optimism to other views regarding the status of CH in a more or less systematic manner. 

\paragraph{7. Skepticism.} We may generally talk about \textit{skepticism} or \textit{pessimism} in order to include views that regard CH either as a legitimate but already solved problem, as a pseudo-problem or ill-posed question and those regarding Gödel's program as unfeasible. Below we deal with what we consider the main arguments for the skeptic position.

\subparagraph{a. CH has already been solved.} Perhaps the first view against what today is vindicated as Gödel's program was given by Cohen himself, regarding the proof of the undecidability of CH in ZFC as the solution to the problem. Indeed, if we take ZFC as providing the standard theory of sets then the undecidability of CH only shows, at best, the limitations already announced by Gödel's incompleteness theorems, namely, that CH has no definite answer. As a consequence, the alleged problem in determining whether, say, ZFC + CH or ZFC +  $\neg$CH is \textit{the} correct system disappears. Additionally, from a more general attitude towards the undecidability phenomenon we may simply argue that we are asking the incorrect questions in set theory and that, perhaps, this explains why such phenomenon has proven to be so pervasive in this mathematical area \cite{shelah}.
\\

\textbf{Observation.} First, even though it is interesting to rephrase problems like CH in other settings such as the pcf-theory from \cite{shelah}, the optimist will of course defend the usual formulation of CH as given. After all, we do not try to rephrase other mathematical conjectures in a different way just to solve their modified form. Second, even without holding any form of realism, that new axioms become natural (be it intrinsically or extrinsically)  may be considered as a hard fact from experience. The most well-known example is that of AC. Therefore, we cannot take as \textit{a priori} granted that ZFC will not be augmented in the future and that, hence, that CH will remain \textit{absolutely undecidable} (see \cite{koellner2}). This would certainly be `an apriorism with its sign reversed' (cf. 2.a.5). Ultimately, while the optimist can recognize that the problem of the decidability of CH in ZFC is indeed solved, this is just the initial condition of the new general problem.
\\

\textbf{Moral.} This position takes the problem as being `decide CH in ZFC' but nothing prevent us from generalizing it into `decide CH in some natural extension of ZFC'. Hence, Gödel's program should be better understood from this general viewpoint.

\subparagraph{b. CH is vague.} Another view is that CH is not a legitimate mathematical problem because it is \textit{vague}. As Feferman argues in \cite{feferman}, CH assumes for its formulation several basic notions (the continuum, subsets of the continuum and mappings between such subsets) that admit differing conceptions, all superficially homogenized by the set-theoretic conception. In particular, the notion of `arbitrary set' is \textit{inherently vague}, meaning that `there is no way to sharpen it without violating what it is supposed to be about' (see \cite{feferman}). Then, CH is void of definite content, it is not even absolutely undecidable. The only way to decide when a given mathematical notion is definite is, according to Feferman, restricting ourselves to weaker systems (e.g. predicative ones). Moreover, the use of increasingly exotic assumptions that appear in several programs for deciding CH leads to problems, both metaphysical (see \cite{fefermannewaxiomsmaddy}) and mathematical since, for instance, `even though the experts in set theory may find such assumptions compelling from their experience of working with them [...], the likelihood of their being accepted by the mathematical community at large is practically nil' (see \cite{feferman}). 
\\

\textbf{Observation.} Feferman is arguing for a position that goes in a directly opposite direction to our phenomenological optimism. Indeed, he believes that we cannot elaborate further properties of the idea of set without making (synthetic) judgments that may very well be otherwise. The difference between our points of view is, then, that we leave the possibility of obtaining new promising analyses open while he rejects this approach altogether. This is merely a clash between two essentially contradictory beliefs. 
\\

On the other hand, Feferman's treatment of vagueness is certainly philosophically fruitful: what is the nature of the phenomenological analysis that we propose? For instance, if the analysis of `set' and `inclusion' are to be finished, how do we properly recognize that they have not been already finished? If, on the contrary, they stand in need of further work, does this undermine the definiteness of such notions? Perhaps the solution consists in recognizing that this open-ended nature has to do with our very mathematical intuition that keeps evolving. Gödel seemed to believe so (cf. 6.a.2). 
\\

\textbf{Moral.} We can acknowledge that CH is indeed vague if we do not treat as it is, as a mathematical statement formulated in a concrete formal system. But we have already seen in \S7.a how we can take the decision of CH to be a general problem, quantifying over potential natural extensions of ZFC. This helps us vindicating some kind of `flavor' that CH enjoys, as do other open mathematical questions. The key difference is that, say, Goldbach's conjecture is expected to open new possibilities \textit{upwards} (i.e. in the direction of new theorems and definitions) while CH is expected to do so \textit{downwards} (i.e. in the direction of new axioms)\footnote{See \cite{yogp1}.}. In other words, it does not seem natural -- in principle -- to think that a \textit{proof} of Goldbach's conjecture in its present state will involve adding new axioms, be it because natural numbers are sharper than sets or not: the undecidability phenomenon leads to a different understanding of axioms and formal systems in general. Still, we may admit that a non-solution to CH has a different evidence degree than its mere negation (say, adding the failed candidates $V=L$ or PD \textit{versus} doing so with CH or $\neg$CH)\footnote{\textit{ibid.}}.

\subparagraph{c. Set theory needs new primitive notions.} It is well known that Freiling's axiom of symmetry (see, e.g., \cite{jdh, hauseraxioms}) was introduced through an interesting mental experiment regarding the intuitions that we have about randomness and the continuum. Notoriously, this axiom can be proved to be equivalent to $\neg$CH in ZFC. More generally, Mumford has claimed that `if we make random variables one of the basic elements of mathematics, it follows that the CH is false and we will get rid of one of the meaningless conundrums of set theory' (see \cite{mumford} and see \cite{hauseraxioms} for a presentation of similar views). The issue with CH, as it is usually conceived as a problem, is that we must search for new axioms that stand as natural properties relating the primitive notions (set and belonging) \textit{exclusively}. The general programmatic approach that Mumford is following consists in adding new axioms for set theory \textit{while enlarging the set of primitive notions}. As he states, CH `is surely similar to the scholastic issue of how many angels can stand on the head of a pin: an issue which disappears if you change your
point of view' (see \cite{mumford}).
\\

\textbf{Observation.} Despite the fact that Freiling's argument has failed to be nowadays recognized as a proper solution to CH (see, e.g. \cite{maddyaxiomsI}), some philosophical problems remain. The point is that the usual arguments against such argument tend to rely on pragmatical and other extrinsic considerations, rather than intrinsic arguments. This, as Hamkins recognizes (see \cite{jdh}) leads to a problem that cannot be eluded: if this happens with Freiling's arguments, which follows the dream solution template, why should it not happen with every other possible candidate? See below (\S7.d) for our view regarding this problematic.
\\

\textbf{Moral.} On the other hand, the openness in our phenomenological conception prevents us from totally eliminating the possibility described by Mumford. That is, we believe it possible that a satisfactory account of the continuum may lead to the addition of new primitive ideas, in the fashion that we have also described in \S7.b. After all, even if we have at our hand the usual construction of the real numbers, we cannot argue \textit{a priori} against other conceptions of the continuum. In this case, however, a further clarification is needed: one should justify how the enlarged set of primitives represents a suitable conception of set as those making use of only the two usual ones. Here, the weight of tradition remains to be lifted and fought but a positive outcome cannot be denied from the outset. At this point, it seems useful to recall the famous Wittgensteinian remark that
\begin{quote}
    7.d.1 It is not single axioms that strike me as obvious, it is a
system in which consequences and premises give one another \textit{mutual} support. \cite[\S 142]{wittcertainty}
\end{quote}

Moreover, Husserl can be also (loosely?) read arguing that
\begin{quote}
    7.d.2 [...] the determinate equivocality can once again be eliminated by \textit{the joint force} [emphasis is mine] of the axioms, so that we are enabled univocally to determine new and ever new elements from given elements (and here that can only mean elements assumed as given and, as it were, named by means of proper names) on the basis of the axioms, and consequently to regard them likewise as given.  \cite[K I 26/43]{husserlDV}\footnote{\textit{ibid.}}
\end{quote}

\subparagraph{d. Gödel's program is unfeasible.} Feferman is attacking the phenomenological part of Gödel's program (see \cite{fefermannewaxiomsmaddy}) and, thus, his view finds the overall approach to CH problematic. Hamkins \cite{jdh} has also worked out an argument pointing out how the phenomenological and mathematical parts of Gödel's program seem to clash. He argues that, nowadays, set theorists have experience working with both CH and $\neg$CH. Hence, if someone claimed to have, say positively, decided CH on the basis of some new, evident axioms, the \textit{hard} experience working with $\neg$CH would undermine such evidence: on the one hand set theorists would allegedly accept such axioms while on the other their experience would prevent them from such approval. In this case, as he puts it, `before we will be able to accept CH as true, we must come to know that our experience of the $\neg$CH worlds was somehow flawed' (see \cite{jdh}).
\\

\textbf{Observation.} The Gödelian could argue against such view on the basis of 2.b.2 in the following way: such parallel experience with both CH and $\neg$CH is a matter of superficiality of the analysis of the primitive notions of the formal systems. From our phenomenologically optimistic standpoint, both experiences are perfectly understandable as corresponding to different analyses of the notion of set. In this way, one analysis will deem CH as (following from) some potential essential feature of such notion and the other will do so with $\neg$CH. Of course, in the concrete case of Gödel's approach, this would mean that (only) one of these analysis would correspond to the iterative conception of set (or, alternatively, his conception of sets as guided by Ackermann's principle). The problem has to do with the meaning of `the intended theory sets'. Pluralism argues that either there is no such a thing or that there are multiple such theories. The monist can, after all, adopt a more modest position and simply try to provide the corresponding analyses of the preferred notion of set. Our phenomenological optimism grants this as a plausible alternative.
\\

In a nutshell, Hamkins seems to argue that Gödel's ideal solution consists in two, temporally differentiated steps. Rather, as Wang points out regarding the vicious circle principle, `this matter of presupposition [...] is not a question of temporal priority' (see \cite[8.2]{wang}). Gödel's \textit{dream solution} (as it is put in \cite{jdh}) involves that the evidence of a new axiom will, immediately, either include the experience working with decided statement or surpassing the one against it (cf. 1.b.3). If we can \textit{perceive} a possibility of solution (see 3.b.2), then such perception is a simple act, phenomenologically speaking\footnote{See \cite{yogp2}.}.
\\

\textbf{Moral.} If the mathematical and the phenomenological are not properly differentiated and we do not concede an equal amount of attention to them or, as Gödel puts it, if we wish to deal with formal systems disregarding the informal notions behind them, we will likely end up confronting equally valid experiences and equally reasonable pragmatic values\footnote{Again, this insight was also presented in \cite{yowitt}.}. But if we appropriately limit the field of our investigations in a legitimate way regarding such notions, then there is no reason to assume that experience can precede evidence.

\paragraph{8. Extrinsic optimism.} Regarding the specific case of CH, we distinguish \textit{extrinsic optimism} from \textit{intrinsic optimism} in that the former intends to settle CH, as a legitimate mathematical problem, through new, extrinsically justified axioms, where the latter also accepts intrinsically sanctioned candidates or, at least, regards the task of providing these as feasible. Our phenomenological optimism is, of course, a form of the latter.

\subparagraph{a. CH is purely mathematical, I.} As we have seen, Feferman centers his arguments around the notion of intrinsic justification. But then, the door of extrinsic criteria remains open. Maddy \cite{maddyaxiomsI, fefermannewaxiomsmaddy} has developed a view in which CH acquires the status of a legitimate mathematical question but in which extrinsic justifications are preferred to the intrinsic ones. Of course, for the Gödelian, the notion of intrinsic naturalness has to do with the very purpose of set theory, namely, to provide a foundation of mathematics. Despite this, Maddy has been able to provide a different account of what a `foundation' is expected to do \cite{maddyfoundations}. Roughly, the core of her set-theoretical naturalism consists in regarding set theory as providing a `unified arena in which set theoretic surrogates for all classical mathematical objects can be found and the classical theorems about these objects can be proved' (see \cite{fefermannewaxiomsmaddy}), which stands far away from either a metaphysical or an epistemological grounding of mathematical knowledge. Hence, Feferman arguments may work for dismantling the intrinsic side of Gödel's quest, but not the extrinsic one. Maddy's rephrasing of Gödel's program is the following: ‘extrinsic reasons are easy to recognize; they involve the consequences of a given axiom candidate, its fruits, if you will’, so ‘we need to assess the prospects of finding a new axiom that is well-suited to the [purely mathematical, methodological] goals of set theory and also settles CH’ (see \cite{fefermannewaxiomsmaddy}).
\\

\textbf{Observation.} First, we believe that establishing the supremacy of either extrinsic or intrinsic criteria leads to a wrong understanding of Gödel's program. Recall, for example, that in \S 3.c we saw how the same axiom candidate may be justified in this two ways simultaneously. It may very well be the case that the notions of intrinsic and extrinsic are two ends of a continuum. In order to vindicate this possibility, consider one of the extrinsic criteria that Maddy proposes: MAXIMIZE (see, e.g., \cite{maddyfoundations}). Gödel famously remarked that, for instance, ‘only a maximum property would seem to harmonize with the concept of set' from 1.a.2. Hence, an axiom such as $V=L$ would be sanctioned as unnatural on this ground. Now, Maddy claims that this exclusion is carried out from a completely methodological point of view while Gödel certainly concluded this from his pre-theoretical conception of sets. The tendency to collapse both kind of criteria is explained as follows: when we take intrinsic naturalness as essential, extrinsic naturalness may be considered as ultimately relying upon it, as if it was its hidden cause; when, on the contrary, extrinsic naturalness is taken to be essential, intrinsic naturalness seems a mere philosophical ornament placed on top of a purely pragmatical issue. This tendency is to be attacked under the title of reductionism. On the other hand, while we have explored intrinsic naturalness in this paper, a similar task remains for the extrinsic case. Similar questions, if not even more difficult, will surely arise in such context\footnote{For instance, one could also ask whether there are other forms of justifications of axioms. C. Straffelini once suggested in private communication that perhaps a new form of naturalness may arise by paying attention to the possible consequences that a new axiom may have in the natural sciences. Even if this could be considered as a concrete type of extrinsic justification, it certainly involves considerations of a different kind.}.
\\

Second, erasing the epistemological component of Gödel's program leads to a collapse between the two halves that compose it, as we have seen in \S 4, and this is untenable from the point of view of both our phenomenological optimism and Gödel's original intentions. But `erasing' is certainly not the appropriate term for describing Maddy's strategy. It seems more correct to claim that she is hiding the epistemological worries regarding set theory under the carpet of the purely methodological and mathematical tasks developed by the working set theorists. From the classical foundationalist view of mathematics, Maddy is taking set theory to be both a \textit{foundation} and a \textit{branch} by including foundational goals of the theory along the merely methodological ones (see \cite{maddyfoundations}). This means that certain philosophical features of set theory are defended without further comment under the tile of common practice which amounts to a form of dogmatism\footnote{For the foundation-branch duality, see \cite{yoch, yogp1}.}. The only possible solution to this impasse would consist in recognizing the proper place of the epistemological (against the merely methodological) in set theory, which seems a quite serious task. What remains clear is, however, that Maddy's desire of extirpating the epistemological from the methodological leads to surreptitiously collapsing intrinsic and extrinsic naturalness.
\\

Third, even if we did center ourselves exclusively in extrinsic criteria, further elucidating what the methodological goals of the discipline are would consist in an epistemological task, which by definition is precluded from set theory, according to Maddy. Hence, we either concede self-evidence to such pragmatic principles (and recall that self-evidence is synonymous with mysticism and dogmatism, cf. \S3.a) or we admit that a phenomenological analysis, similarly to the one that Gödel proposes for the primitive notions, is to be carried out for clarifying these methodological principles, which stands as a reasonable and attractive project by itself\footnote{For the problem of mysticism, see \cite{yodfl}.}.
\\

Lastly, if dealing with the purely methodological amounts to excluding from our attention the informal insights that we have vindicated in \S7.d then Hamkins' arguments applies with full force: indeed, if we only took experience as relevant then the resulting situation confronting two experiences would not be susceptible of a satisfactory decision and, thus, Maddy's ideal solution to CH would hardly be possible.
\\

\textbf{Moral.} Precluding intrinsic naturalness from Gödel's program arises several problems: first, eliminating several natural examples; second, collapsing intrinsic and extrinsic naturalness; third, conceding certain methodological principles the status of \textit{a priori} values and lastly, enabling Hamkins' argument above. The alternative is to maintain \textit{both} intrinsic and extrinsic justifications as providing a field of possibilities, not of exclusions. Recognizing the proper place of these justifications in the overall program, when carried out specifically, belongs to the overall phenomenological task that we argue for from our phenomenological optimism.

\subparagraph{b. CH is purely mathematical, II.} As it is widely known, one of the most serious attempts for settling CH -- together with Foreman's program -- has been provided by Woodin (see \cite{bellotti, rittberg, woodinch2}). In few words, Wooding has attempted two strategies regarding CH: (a) providing grounds for $\neg$CH through the $\Omega$-conjecture and (b) providing grounds for CH through structural properties of the so-called Ultimate-$L$ program. With a similar spirit to that of Maddy, Woodin's approach depends on his \textit{mathematical traction}, i.e. his belief that (at least some) philosophical questions can be settled by mathematical means (see \cite{rittberg}). This strategy seems more appealing and refreshing than Gödel's program, which remains entangled with \textit{demodé} epistemological questions (see \cite{fefermannewaxiomsmaddy}), since it connects purely mathematical properties of other theories to the solution of CH. One positive property that this approach shows to have is that, even if it has not settled CH \textit{yet}, it at least shows `convincing evidence that there is a solution' \cite{woodinch2}. Woodin also argues against the skeptic or, what he calls, the `widget possibility', that is, `the future where it is discovered that instead of sets we should be studying widgets' and in which `the axioms for widgets are obvious and, moreover, that these axioms resolve the Continuum Hypothesis (and everything else)' (see \cite{woodinch2}).
\\

\textbf{Observation.} First, Gödel's phenomenological program can be translated for Woodin's approach: intrinsic justifications should be provided for, say, accepting certain features of the theoretical frameworks in which Woodin develops that possible solution of CH. In other words, Gödel's program wishes to deal with some degree of \textit{evidence} and, as such, it can be asked how Woodin's alternative theories are related to that same evidence\footnote{See \cite{yogp1}.}. This is also a concrete application of the preceding arguments in \S8.a.
\\

Second, Woodin is known to profess some sympathy for the view that philosophy provides certain \textit{inspiration} for mathematics and, in particular, set theory\footnote{During a panel discussion in \textit{What can Philosophy do for Set Theory?}, in Barcelona, March 2023.}. In contrast with this positive relation, mathematical traction is a purely negative one. The main difficulty with this is that, as we have seen in \S8.a, one should first clarify what we understand by `philosophical'. Our demarcation between phenomenological and mathematical discourse (cf. \S5.c) certainly clashes with Woodin's conflicting link between both. Moreover, it remains the question on how Woodin's program stands with respect to others seeking a solution to CH, and the previous comments on pragmatical values from \S8.a also apply here.
\\

Finally, Woodin's identification of the `widget possibility' with skepticism is strange from the point of view of phenomenological optimism. Indeed, it clashes with the openness of the phenomenological method and how we should be open to a possible situation in which sets, as given from a concrete conception, rely upon additional primitive ideas\footnote{See \cite{yogp1}.}. Of course, such possibility is not a matter of mere desire (see \S7.c).  
\\

\textbf{Moral.} It is interesting to notice the refinement that the notions of `problem', `solution' and `possibility' have undergone in recent times (see, e.g., \cite{dehornoy}, where a solution for CH is defined in a purely mathematical way). Nevertheless, it seems reasonable to argue that the solution of CH is often confused with its possibility of solution, and these two notions should be understood differently. Woodin, in the quote above, seems to share this symptom\footnote{\textit{ibid.}}. Gödel's phenomenological program deals with the latter notion, while the mathematical side does so with the former. Therefore, similar concerns from \S8.a can be replicated here.

\subparagraph{c. Intrinsic pessimism.} Let us conclude this section with some last observations. Regarding 2.a.5, Wang confesses that he is `not aware of any conspicuous successful examples of definite axioms arrived at' following the phenomenological method (see \cite[5.1]{wang}). Nevertheless, as we have seen in \S3.b, he then provides some philosophical justification of certain axioms of set theory, which could be rephrased, with some work, to constitute intentional analyses of a concrete notion of set. Similarly, what we have seen with Maddy's views in this section suggests that extrinsic optimism may be seen as an optimism regarding the status of CH and, simultaneously, an intrinsic pessimism. 
\\

Hence, this view runs contrary to our phenomenological optimism. It is surprising, however, that both the skeptic and the (usual) extrinsic optimist coincide in arguing against intrinsic naturalness. The former, because she sees something essentially wrong in it and, at the same time, considers it to be central to the overall optimist picture; the latter, because she disregards it as problematic and unessential. From this perspective, phenomenological optimism has more ties to this skeptic than to the extrinsic optimist.

\paragraph{9. Intrinsic optimism.} Lastly, we compare our views to other important approaches belonging to the intrinsic optimistic view.

\subparagraph{a. Vindication of phenomenology.} As far as we know, Hauser's work \cite{hauserch, hauseraxioms, hauserchoice, hauser, hauserintuition} is remarkable for being perhaps the only one concerned with a phenomenological examination of the notion of intrinsic naturalness. In \cite{hauser} he presents an illuminating outline of Gödel's thought regarding CH, where the distinction between ontological and epistemological realism that we have followed above (cf. \S2.b) is made. In \cite{hauserchoice} and \cite{hauserintuition}, some Husserlian remarks are put into practice for determining whether AC is self-evident and for dealing with the problematic of mathematical intuition, respectively. Lastly, \cite{hauserch, hauseraxioms} are directed against both the Feferman's views sketched above (cf. \S7.b) and against the status Freiling's argument as a proper proof of CH and other related views (cf. \S7.c). 
\\

As the reader can infer from this, Hauser's views and ours will differ in few points. One of these is that Hauser wishes to defend a form of epistemological realism: he claims that Gödel is able to avoid the metaphysical issues of mathematical realism by means of a phenomenological account of mathematical intuition (see, e.g., \cite{hauser, hauseraxioms}). From our weaker point of view, such claim remains in need of further elucidation. Even so, as we have argued, epistemological realism requires additional components to a mere optimism regarding the phenomenological method (cf. \S6.c).
\\

Another point of disagreement is explained by the openness that, as we argue, Gödel's phenomenological program should bear. More precisely, regarding Mumford's general proposal of adding new primitive ideas (cf. \S7.c), Hauser argues for the adequateness of ZFC as a `formal framework of mathematics' (see \cite{hauseraxioms}) on two grounds: that every mathematical statement is expressible in the language of ZFC, i.e. \textit{exclusively} using set and belonging as primitive notions; and that every classical theorem is derivable in ZFC. The problem with this line of argumentation is that it relies on essentially extrinsic, i.e. methodological, facts. Certainly, the two asserted reasons are not mathematical facts susceptible of proof; their status is probably similar to that of CTT. But they do not constitute a deep analyses of the notions of set theory: what is to be expected instead is the task that we have described in \S7.c but in the opposite direction\footnote{\textit{ibid.}}. 
\\

\subparagraph{b. Reflection principles.} In \cite{bagaria} one can find perhaps the most systematic exposition of Gödel's criteria for the addition of new axioms of set theory (apart from \cite{wang} and \cite{hauser}). The purpose of this paper is to first present a satisfactory collection of `meta-axioms of set theory' that allow to decide the naturalness of a given candidate and then study particular examples. In particular, it is claimed that bounded forcing axioms are indeed natural axioms of set theory and, hence, that $\mathfrak{c} = \aleph_2$ is a natural solution to CH.  
\\

But here we are mainly concerned with the philosophical details in the aforementioned exposition, for these will shed light on the status of the purported solution. The approach of \cite{bagaria, bagariareflection} coincides in many ways with the one presented in \cite{hauserabsolute}, namely, continuing the task left open by Gödel regarding the refinement of the iterative conception of set by means of Ackermann's principle. In other words, the main intention is to justify the regulative principles concerning the addition of new axioms through the Cantorian conception of the Absolute \cite{jané} which, as we have already seen in \S3.c, implies that reflection principles gather special attention (cf. 3.c.2) (also, see \cite{koellner, koellner2}). The two main potential problems with this view (already anticipated in \S3.c) are (a) justifying that the Absolute conception leads \textit{univocally} to the formalized reflection principles and (b) justifying the Absolute conception itself. Moreover, if these questions do not receive any answer then, in some sense, Hamkins argument against the quest for achieving \textit{evidence} of new axioms applies to the more philosophically entrenched one regarding satisfactory regulative principles as those presented in \cite{bagaria}\footnote{See \cite{yoch}.}.

\newpage

\section*{IV. Final thoughts}

\paragraph{10. Further problems.} We have presented phenomenological optimism as a view that is strong enough to support some form of optimism about the status of CH. We have also seen how taking the phenomenological method for granted without further inspection, regarding it as superfluous or rejecting it altogether amount to the same thing, namely, intrinsic pessimism. However, some problems that need to be addressed further remain. Here, we will briefly deal with some of the issues that arise against the plausibility of the phenomenological method.

\subparagraph{a. Sudden illumination, revisited.} As we have seen in \S6.b, sudden illumination is a key factor of Gödel's overall optimism. This idea is closely tied with some considerations that we have seen in the preceding section. Namely, it seems that Gödel's template for a solution of CH relies upon some kind of spontaneous \textit{feeling} of the naturalness of some candidate for new axiom. More precisely, it seems that the approach that Gödel defends consists in bringing some deep, unconscious properties of sets to the light of consciousness. This would certainly go against the account of phenomenology from 3.a.1, i.e. as a conscious and systematic method, for it is unclear how a series of organized and structured conscious acts could make the unconscious explicit in a non-spontaneous way. In fact, Wang argues for his intrinsic pessimism by claiming that `the role of the unconscious presents a more fundamental obstacle' to the phenomenological method `than the occasional failure to fully communicate our conscious thought' (see \cite[5.1]{wang}). Not only this, as we have seen before (cf. \S7.d), such psychological experience would be undermined by how (some) set theorists work today with CH and $\neg$CH. 
\\

Against this, the point is not that new axioms will impose themselves in a similar way to physical properties (as, e.g. one can read in \cite{maddy}), so that this objective self-evidence will simply stand above every merely surjective content, since this would be tantamount to a positive transliteration of Hamkins view. Rather, the point is to phenomenologically elucidate how a \textit{natural} axiom could be given for us, that is, in an anti-reductionistic, immediate, descriptive and open way. We \textit{do} relate intentionally with axioms in this way, as the usual examples show. Our phenomenological optimism grants that we will be able to gradually provide a systematic account of such intentional relation, which includes the subjacent ones corresponding to the primitive ideas of set theory, for instance. What we acknowledge in each case is that such relations \textit{are} a thing and, additionally, that we will be able to systematically clarify them. 
\\

Hence, in case we arrived at a natural axiom that decides CH, this would not constitute a psychological or metaphysical fact but, rather, a mere fulfillment of an intentional relation and, hence, an indivisible (ideal) act of recognizing such axiom as natural. In particular, the subconscious conception of sudden illumination regarding the solution of CH should be rejected as directly opposing the methodological values of  the phenomenological method (see, e.g., \cite{scherer}). Phenomenology does not provide new data, it simply takes our knowledge and tries to organize it in a systematic manner. Our optimism tells us that, with prolonged efforts, new relations will be made evident by means of these analyses. 

\subparagraph{b. Communicability} The other main problem that Wang sees in the pursue of the phenomenological method is that of the communicability of results and procedures. Indeed, Gödel's conception of sudden illumination also implies some degree of \textit{private} experience and introspection that would be probably useless when made explicit. It is not surprising, then, that `[n]either the axiom of choice, the axiom of replacement, the ``axiom'' of constructibility, the ``axiom'' of determinacy, nor even Dedekind's axioms for arithmetic were obtained by going back to the ultimate acts and contents of our consciousness in the manner recommended by phenomenology' (see \cite[5.1]{wang}). We cannot even tell if Gödel \textit{made use} of the phenomenological method when he arrived at a Ackermann's principle as a refinement of the iterative conception of set (cf. \S3.d). 
\\

But these seemingly hostile views against phenomenology are merely the usual worries associated with mathematical creativity and creativity in general. The point is not to follow exclusively the phenomenological method (as incarnated in several methodological values) but to systematize our knowledge according to it. This, of course, may be an activity that lies far away from the natural sciences and mathematics but, ultimately, phenomenology has to do with \textit{public} contents. When Wang talks about overview, idealization and extensionalization (cf. \S3.b) he is describing some intentional relations that we experience when dealing with sets (a Kantian would add: when we build them in intuition). That is, we are not looking at the private psychological experiences that are accidental to such relations (such as the use of \textit{concrete} mental images or diagrams). This is where the phenomenological and the grammatical are ultimately connected. If we deal only with these intentions -- which certainly are \textit{not} psychological in nature -- then it is clear how the communicability problem may be overcome.

\paragraph{11. Conclusion.} In this paper we have presented Gödel's program in his original formulation. In section I we have seen how Gödel can be read in several places as defending epistemological realism (every mathematical question can receive a yes-no answer) and, as a consequence, what we have called phenomenological optimism (primitive notions of set theory admit analyses in the course of which new axioms will become evident). In section II we have defended: (i) that the phenomenological method should be understood as the main approach for dealing with the epistemological problems arising from the study of intrinsic criteria and that this was probably what Gödel also believed, (ii) that the epistemological part of Gödel's program should be separated radically from the mathematical one, (iii) that phenomenological investigations should follow certain methodological values and, notoriously, the openness motto, (iv) that Gödel's phenomenological program is to be understood as a general approach to the undecidability phenomenon, (v) phenomenological optimism as our preferred form of intrinsic optimism, which is strictly weaker than epistemological objectivism. Finally, in section III we have contrasted our view with alternative forms of intrinsic optimism, as well as others that we have labeled as pessimistic and as extrinsically optimistic.

\begin{comment}

\newpage

\textcolor{red}{Wang sobre la fen: phenomenological method to discover the axioms for the primitive
concepts of philosophy (and of more restricted fields). But I am not
aware of any conspicuous successful examples of definite axioms arrived
at in this manner. Neither the axiom of choice, the axiom of replacement,
the "axiom" of constructibility, the "axiom" of determinacy, nor even
Dedekind's axioms for arithmetic were obtained by going back to the
ultimate acts and contents of our consciousness in the manner recommended
by phenomenology (5.1)}

\textcolor{blue}{Hauser sobre AC. Sobre maddy y a veces hauser: 
We have argued, therefore, that Gödel's phenomenological program connects more naturally with certain form of epistemological objectivism (2.b.4) than with an ontological or strong version of such position (1.a.1). This is because in the phenomenological method, not only positive facts regarding natural sciences (and here we may include psychology) are precluded but also those belonging to metaphysics and ontology. }
posibilidad de criterios empíricos (no extrínsecos...) concetar con CH y la CTT.

\end{comment}

\newpage
%usas \ref{quant} donde quant es la label de una entrada de la bibliografía

%\bibliographystyle{abbrv}
%\bibliographystyle{plainnat}

\end{document}